\documentclass[11pt,a4paper]{amsart} 

\usepackage{latexsym}

\usepackage[a4paper , lmargin = {2.5cm} , rmargin = {2.5cm} , tmargin = {2.5cm} , bmargin = {2.5cm} ]{geometry} 

\usepackage{geometry}
\usepackage{pinlabel}         
\usepackage{graphicx}
\usepackage{amssymb, amsthm}
\usepackage{epstopdf}
\usepackage{romannum}
\usepackage{tikz}
\usepackage{subfig}
\usepackage{MnSymbol}
\usepackage[arrow, matrix, curve]{xy}

  	\newcommand{\Z}{\ensuremath{\mathbb{Z}}}

   	\def\Aut{{\rm{Aut}}$(F_n)$}
        \def\Autn{{\rm{Aut}}(F_n)}
        
   	\def\SAut{{\rm{SAut}}$(F_n)$}
        \def\SAutn{{\rm{SAut}}(F_n)}
	\def\GL{{\rm{GL}}}

	\def\Z{{\mathbb{Z}}}

\theoremstyle{plain}
\newtheorem*{NewTheorem}{Theorem}

\newtheorem{theorem}{Theorem}[section]

\newtheorem{proposition}[theorem]{Proposition}

\newtheorem{remark}[theorem]{Remark}

\newenvironment{pf}{\par\medskip\noindent\textit{Proof.}~}{\hfill $\square$\par\medskip}

\title[\SAut\ cannot act on small tori]{$\mathbf{SAut}\boldsymbol{(F_n)}$ cannot act on small tori}
\author{Olga Varghese}
\date{\today}
\address{Olga Varghese\\
Department of Mathematics\\
M\"unster University\\ 
Einsteinstra\ss e 62\\
48149 M\"unster (Germany)}
\email{olga.varghese@uni-muenster.de}

\begin{document}
\pagenumbering{arabic}
\begin{abstract}
We study smooth actions of \SAut, the unique subgroup of index two in the auto-\\morphism group of a free group of rank $n$,  as a part of the generalized 'Zimmer program'. In particular, we show that every action of~\SAut\ on a low dimensional torus is trivial.
\end{abstract}

\maketitle

\section{Introduction}
In the mathematical world, this work is located in the area of geometric group theory, a field at the intersection of algebra, geometry and topology. Geometric group theory studies the interaction between algebraic and geometric properties of groups. One~is interested to understand on which 'nice' geometric spaces a given  group can act in a reasonable way and how geometric properties of these spaces are reflected in the algebraic structure of the group.
Here, the spaces will  be tori,  while the group will be \SAut. 

In this work we study smooth actions of \SAut\ as a part of the generalized 'Zimmer program' whose aim is to understand actions of large groups on compact manifolds, see e.g. \cite{Zimmer}.  Given a group whose low dimensional linear representations are trivial, one wants to know whether all group representations into~${\rm Diff}(M)$, where $M$ is a low dimensional compact manifold, are also trivial. 

To be precise, let $\mathbb{Z}^n$ be the free abelian group and $F_{n}$ the free group of rank $n$. One goal for a group theorist is to unterstand the structure of their automorphism groups, 
${\rm GL}_n(\mathbb{Z})$ resp. \Aut.
The abelianization map $F_{n}\twoheadrightarrow \Z^{n}$ gives a natural epimorphism 
$\Autn\twoheadrightarrow\GL_{n}(\Z).$
The special automorphism group of $F_n$, which we will denote by \SAut, is defined as the preimage of ${\rm SL}_{n}(\Z)$ under this map. 
Much of the work on \Aut\ and \SAut\ is motivated by the idea that ${\rm GL}_n(\mathbb{Z})$ and \Aut\ resp. ${\rm SL}_n(\Z)$ and \SAut\  should have many properties in common. 
Here we follow this idea and present analogies between these groups with respect to smooth actions on low dimensional tori.

It was shown by \textsc{Weinberger} in \cite{Weinberger} that any smooth action of ${\rm SL}_n(\mathbb{Z})$ on the $d$-dimensional torus $\mathbb{T}^d=\mathbb{S}^1\times\ldots\times\mathbb{S}^1$ is trivial if $d< n$. This bound on the dimension of the torus is sharp, as
${\rm SL}_n(\mathbb{Z})$ admits a linear faithful action on the flat torus $\mathbb{R}^n/\mathbb{Z}^n$.
We prove an analogue of \textsc{Weinberger's} result for smooth actions of \SAut.

\begin{NewTheorem}
Let $n\geq 3$ and 
$\Phi:\SAutn\rightarrow{\rm Diff}(\mathbb{T}^d)$
be a smooth action. If~$d<n$, then $\Phi$ is trivial.
\end{NewTheorem}
Combining one result of \textsc{Weinberger} concerning smooth actions of non-abelian finite groups on tori with results of \textsc{Bridson} and \textsc{Vogtmann} about strong constrains on homomorphisms from \SAut, we prove the above theorem.

\section{The automorphism group of a free group}
As the main protagonist in this work is the group \SAut. We start with the definition of this group, and establish some notation to be used throughout. 

We begin with the definition of the automorphism group of the free group of rank~$n$.
Let $F_{n}$ be the free group of rank $n$ with a fixed basis $X:=\left\{x_{1}, \ldots, x_{n}\right\}$. We denote by \Aut\ the automorphism group of $F_n$ and by \SAut\ the unique subgroup of index two in \Aut. 

Let us first introduce a notations for some elements of \Aut.   We define  involutions $(x_{i},x_{j})$ for $i, j=1,\ldots, n$,~$i\neq j$ as follows:
\[
\begin{matrix}
(x_{i},x_{j})(x_{k}):=\begin{cases} x_{j} & \mbox{if $k=i$,}  \\ 
			     x_{i} & \mbox{if $k=j$,}  \\ 
			     x_{k} & \mbox{if $k\neq i,$ $j$.}  
\end{cases}
\end{matrix}
\]
In particular, the alternating group ${\rm Alt}_n$ is a subgroup of \SAut.

The following variant of a result by \textsc{Bridson} and \textsc{Vogtmann} \cite[3.1]{Vogtmann} will be used here to prove that certain actions on spaces of \SAut\ are indeed trivial. For a detailed proof the reader is referred to \cite[1.13]{Diplomarbeit}.
\begin{proposition}
\label{faktorisiertSL}
Let $n\geq 3$, $G$ be a group and $\phi : \SAutn\rightarrow G$ a group homomorphism. If there exists $\alpha\in{\rm Alt}_n-\left\{{\rm id}_{F_n}\right\}$ with $\phi(\alpha)=1$, then $\phi$ is trivial.
\end{proposition}

It was proven by \textsc{Bridson} and \textsc{Vogtmann} in \cite[1.1]{Vogtmann} that if $n\geq 3$ and $d<n$, then \SAut\ cannot act non-trivially by homeomorphisms on any contractible manifold of dimension $d$. Hence 
\begin{proposition}
\label{Dipl}
Let $n\geq 3$ and $\rho:\SAutn\rightarrow{\rm SL}_{d}(\mathbb{R})$
be a linear representation. If $d<n$, then $\rho$ is trivial.
\end{proposition}
In \cite{Varghese} we proved in a purely group theoretical way the above proposition.

We note that the group \SAut\ is perfect, therefore the image of a linear representation of \SAut\ is a subgroup of ${\rm SL}_{d}(\mathbb{R})$.

\section{Proof of main theorem}
We begin the proof of main theorem by the following proposition by \textsc{Weinberger}. The proof of this result is decidedly smooth. First, we note that the first singular cohomology group of a torus ${\rm H}^1(\mathbb{T}^d;\mathbb{Z})$ is isomorphic to $\mathbb{Z}^d$. 
\begin{proposition}(\cite[2]{Weinberger})
\label{ToriAb}
Let $G$ be a non-abelian finite group and $\psi: G\rightarrow{\rm Diff}(\mathbb{T}^d)$ a smooth action. If 
\begin{align*}
{\rm H}^1(\psi): G&\rightarrow{\rm GL}_d(\mathbb{Z}),\\
g&\mapsto {\rm H}^1(\psi(g)): {\rm H}^1(\mathbb{T}^d;\mathbb{Z})\rightarrow {\rm H}^1(\mathbb{T}^d;\mathbb{Z})
\end{align*} 
is trivial, then $\psi$ is not injective.
\end{proposition}

Now we have all the ingredients to prove  
\begin{theorem}
Let $n\geq 3$ and 
$\Phi:\SAutn\rightarrow{\rm Diff}(\mathbb{T}^d)$
be a smooth action. If~$d<n$, then $\Phi$ is trivial.
\end{theorem}
\begin{pf}
According to Proposition \ref{Dipl}, the following action is trivial for $d<n$:
\begin{align*}
\iota\circ{\rm H}^1(\Phi):\SAutn&\rightarrow{\rm SL}_d(\mathbb{Z})\hookrightarrow{\rm SL}_d(\mathbb{R})\\
\alpha&\mapsto{\rm H}^1(\Phi(\alpha))\mapsto{\rm H}^1(\Phi(\alpha)).
\end{align*}
In particular, the map ${\rm H}^1(\Phi)$ is trivial.

We consider the subgroup ${\rm Alt}_n$ in \SAut\ and the restriction of $\Phi$ to this group. The map  ${\rm H}^{1}(\Phi_{|{\rm Alt}_n})$ is trivial, therefore by Proposition \ref{ToriAb} the map $\Phi_{|{\rm Alt}_n}$ is not injective and by Proposition~\ref{faktorisiertSL} it follows that $\Phi$ is trivial.
\end{pf}

\begin{remark}
Note that this bound on the dimension of the torus is sharp, as
${\rm SL}_n(\mathbb{Z})$ admits a linear faithful action on the flat torus~$\mathbb{R}^n/\mathbb{Z}^n$.
Therefore \SAut\ admits a smooth non-trivial action on the flat torus as well:
\[
\SAutn\twoheadrightarrow{\rm SL}_n(\mathbb{Z})\rightarrow{\rm Diff} (\mathbb{R}^n/\mathbb{Z}^n).
\]
\end{remark}

\end{document}